\title{The Complexity of Pebbling and Cover Pebbling}
\author{
Nathaniel G.~Watson\\
Washington University in St.~Louis\\
\tt{ngwatson@wustl.edu} }
\date{April 20, 2005}
\documentclass [12 pt]{article}
\usepackage{amsfonts}
\usepackage[dvips]{graphicx}
\usepackage{amsthm}
\usepackage{amsmath}
\usepackage{amssymb}
\begin{document}
\def\diam{{\rm diam}}
\def\ess{\rm ess}
\def\p{{\mathbb P}}
\def\ep{\varepsilon}
\def\P{{\rm Po}}
\def\cf{{\cal F}}
\def\cl{{\cal L}}
\def\e{{\mathbb E}}
\def\v{{\mathbb V}}
\def\l{\lambda}
\def\dist{{\rm dist}}
\def\lr{\left(}
\def\rr{\right)}
\def\cd{\cdot}
\def\ts{\thinspace}
\def\lc{\left\{}
\def\rc{\right\}}
\def\qed{\vbox{\hrule\hbox{\vrule\kern3pt\vbox{\kern6pt}\kern3pt\vrule}\hrule}}
\newcommand{\hyp}{\mathcal{H}}
\newcommand{\hypp}{\mathcal{K}}
\newcommand{\Nu}{$\begin{Large}$\nu$\end{Large}$}
\newcommand{\ignore}[1]{}

\newtheorem{thm}{Theorem}
\newtheorem{defn}[thm]{Definition}
\newtheorem{result}{Result}
\newtheorem{lemma}[thm]{Lemma}
\newtheorem{cor}[thm]{Corollary}
\newtheorem{prop}[thm]{Proposition}
\maketitle

\begin{abstract}This paper discusses the complexity of graph pebbling, dealing
with both traditional pebbling and the recently introduced game of
cover pebbling. Determining whether a configuration is solvable
according to either the traditional definition or the cover pebbling
definition is shown to be $N\!P$-complete. The problem of
determining the cover pebbling number for an arbitrary demand
configuration is shown to be $N\!P$-hard. \footnote{The significant
results of this paper, with the exception of Theorem 11, have been
independently and simultaneously obtained by K. Milans and B. Clark
in \cite{Milans}, which also analyzes the complexity of optimal
pebbling and of determining the pebbling number of a graph.}
\end{abstract}

\section{Introduction}

Graph pebbling, first suggested by Lagarias and Saks, has been the
subject of many recent developments. It was first introduced into
the literature by Chung in \cite{Chung}, and has been developed by
many others including Hurlbert, who published a survey of pebbling
results in \cite{Glenn2}.

Given a graph $G$ (we assume all graphs we consider in this paper
are connected), we imagine that we can distribute pebbles on its
vertices in arbitrary arrangements called \emph{configurations}.
Formally, a configuration $C$ on a graph $G$ is a function
$C:V(G)\rightarrow \mathbb{N}\bigcup \{0\} $ which represents such a
distribution. We say a configuration $C$ on a graph $G$
\emph{contains} a configuration $C'$ if $C(v)\geq C'(v)$ for all $v
\in G.$ Define the $\emph{size}$ $|C|$ of a configuration $C$ to be
the total number of pebbles, that is, $\sum_{v \in G} C(v).$ By a
\emph{pebbling move}, we will mean a change made to a configuration
of pebbles by taking two pebbles from some vertex which has at least
two pebbles and placing one pebble on an adjacent vertex. We call a
vertex \emph{reachable} from $C$ if there is some sequence of
pebbling moves starting at $C$, which ends in a configuration which
has at least one pebble on this vertex.

The most commonly studied pebbling questions are about a concept
called the ``solvability'' of a configuration. We first consider the
traditional concepts, which we will distinguish with the adjective
``canonical" to avoid confusion. Call a configuration $C$ on a graph
$G$ \emph{canonical pebbling solvable} if every vertex in $G$ is
reachable from $C$. Define the \emph{canonical pebbling number},
$\pi(G)$ to be the smallest integer $k$ such that every
configuration of size $k$ on $G$ is canonical pebbling solvable.

This paper will focus primarily on a concept called ``cover
pebbling" which was introduced by the authors of \cite{Glenn1}, but
suggested earlier by Herscovici in \cite{Herv}. Given a graph $G$,
we imagine not only a configuration $C$ but a \emph{demand}
configuration $D$ on $G.$ We say $C$ is \emph{cover solvable} or
simply \emph{solvable} for $D$ if it is possible, through a sequence
of pebbling moves, to move to from $C$ to a configuration $C'$ which
contains $D.$ We define the \emph{cover pebbling number}
$\gamma_G(D)$ to be the smallest integer $k$ such that ever
configuration of size $k$ on $G$ is cover solvable for
$D.$\footnote{We break from the terminology in \cite{Glenn1} which
calls a demand configuration a ``weight'' configuration, denoted by
``$w$''. Also, we use a different notation for the cover pebbling
number by reversing the role of the subscript and the argument.}

Denote by $U$ the configuration for which $U(v)=1$ for all $v \in
G,$ which we will call the \emph{unit} configuration on $G$ and
denote by $R_v$ the configuration for which $R_v(v)=1$ and
$R_v(w)=0$ for all $w \not= v,$ which we will call the
\emph{$v$-reachability} configuration, because clearly $C$ is
solvable for $R_v$ only if $v$ is reachable from $C.$ Note that
therefore a configuration $C$ is canonical pebbling solvable if and
only $C$ is solvable for $R_v$ for all $v \in G.$

In this paper, we show several results about the complexity of cover
pebbling, including the fact that the question of whether a given
configuration on a graph is solvable for the unit configuration is
$N\!P$-complete. For this, it will be necessary to show that certain
constructed configurations on possibly large graphs are not
solvable. Our first task is to develop tools for proving that a
configuration is not solvable, which will involve slightly
generalizing the concept of solvability.

\section{Preliminaries}

Cover pebbling literature has thus far worked with the definitions
given in the introduction, thinking of weight functions and
configurations as non-negative. However, in \cite{monthlypaper} we
see that if we create a concept called ``negative pebbling" in which
we allow configurations to have any integer number of pebbles on
each vertex, then a non-negative configuration is ``negative
pebbling solvable" (solvable by a sequence of pebbling moves,
allowing negative numbers of pebbles for a given demand
configuration) for a non-negative demand configuration if and only
if it is cover solvable for this configuration as defined in the
introduction, that is, solvable by a sequence of moves with only
non-negative intermediate configurations. Actually,
\cite{monthlypaper} shows this result only for unit configurations,
but the proof applies to arbitrary demand configurations
\emph{mutatis mutandis}, and we repeat this proof for the general
case below.

\begin{thm} \cite{monthlypaper}
Let $G$ be a graph, $C$ a configuration on $G,$ and $D$ a
non-negative demand configuration. Let $|G|=m$ and label the
vertices of $G$ as $v_1,\ldots, v_m$. Then $C$ is cover solvable if
and only if there exist integers $n_{ij}\geq 0$ with $1\leq i,j \leq
m$ and $n_{ij}=0$ and $n_{ji}=0$ whenever $\{v_i,v_j\} \notin E(G) $
such that for all $1\leq k \leq m,$
$$C(v_k)+\sum_{l=1}^m n_{lk}-2\sum_{l=1}^m n_{kl} \geq D(v_k).$$
\end{thm}

\begin{proof}
First, suppose $C$ is solvable for $D$. Then find some sequence of
pebbling moves which solves $C$ for $D$, and let $n_{ij}$ be the
total number of pebbling moves from $v_i$ to $v_j$ in this sequence.
Then after all the moves, there are exactly $$C(v_k)+\sum_{l=1}^m
n_{lk}-2\sum_{l=1}^m n_{kl}$$ pebbles left on $v_k$,
  which is always at least $D(v_k)$ because of the fact that this sequence of moves solves
  $C$ for $D$.

Conversely, suppose such numbers $n_{ij}$ exist. This means that
there {\it does} exist a sequence of moves that solves $C$, with
$n_{ij}$ moves being made from $v_i$ to $v_j$, but possibly with
some illegal ``negative pebbling" along the way. To be explicit, we
can make $n_{12}$ moves from $v_1$ to $v_2,$ then $n_{21}$ moves
from $v_2$ to $v_1,$ then $n_{13}$ moves from $v_1$ to $v_3$ and so
forth, making in general $n_{ij}$ moves from $v_i$ to $v_j$ for all
$1 \leq i,j \leq n$ with $i \not=j$. After the moves are done (which
might require using ``negative pebbles,") we calculate that the
resulting configuration, which we call $\tilde{C},$ is
$$\tilde{C}(v_k)=C(v_k)+\sum_{l=1}^m n_{lk}-2\sum_{l=1}^m n_{kl}\geq
D(v_k)$$ for each $k,$ so $\tilde{C}$ will contain $D.$

We show, however, that it is possible to {\it legally} make enough
moves from our list of $n_{ij}$ moves from $v_i$ to $v_j$ for each
$i,j$ to cover solve the graph. In fact, we may proceed haphazardly,
making moves from our list from a vertex $v_{i'}$ to a vertex
$v_{j'}$ such that less than $n_{i'j'}$ moves from $v_{i'}$ to
$v_{j'}$ have already been made and there are at least two pebbles
on $v_{i'},$ as long as such a pair of vertices $v_{i'}$ and
$v_{j'}$ exists. If no such pair $\{v_{i'}, v_{j'}\}$ is left, then
for each $(i,j)$, either $n_{ij}$ moves have been made from $i$ to
$j$ or else there is at most 1 pebble on vertex $v_i$. Let $C'$ be
the configuration left on $G$ after these moves and $S$ be the set
of $v_i \in G$  for which the total number of moves from $v_i$
having already been made is less than $\sum_{l=1}^m n_{il}.$

If $S=\emptyset$ then clearly for every $1\leq i,j\leq m$ we have
made $n_{ij}$ moves from $v_i$ to $v_j$ and thus, for every $k$
there are  $\tilde{C}(v_k)$ pebbles on $v_k,$ so $C'=\tilde{C}$ and
since $\tilde{C}$ contains $D,$ we have solved $C$ for $D$ through
our sequence of moves.

If  $S \not= \emptyset$ then call the moves that remain to be made
the \emph{inexecutable moves} (clearly all originate from $S$), and
consider the total number of them ($\sum_{i\in S} \sum_{l=1}^m
n_{il}$ minus the number of moves that have already been made from
vertices in $S.$) These moves, if they could executed with negative
pebbling, would transform $C'$ to $\tilde{C}.$ By the definition of
$S,$ there are at least $|S|,$ at least one inexecutable move for
each vertex in $S.$ The total is exactly $|S|$ only if there is
exactly one inexecutable move from each vertex. Also, we know
$C'(v)\leq 1$ for all $v \in S,$ for a total of at most $|S|$
pebbles. Since these moves originate from vertices of $S,$ if
executed, they would each would remove one pebble from $S$ if they
end at a vertex in $S,$ and two if they end at a vertex outside of
$S.$

Thus each inexecutable move must both begin and end in $S,$ for
otherwise $S$ would be left with a negative total number of pebbles
at the end of the pebbling sequence in configuration $\tilde{C}$ ,
which is impossible since  $\tilde{C}$ contains non-negative
configuration $D$ and so is non-negative on all vertices. Even if
all moves begin and end in $S$, however, we end up with at most $0$
pebbles on $S,$ and so must have exactly 0 pebbles on each vertex in
$S,$ that is $\tilde{C}(s)=0$ for all $s \in S.$ Further, since all
inexecutable moves begin and end in $S,$ we know $\tilde{C}$ is
equal to $C'$ on all vertices not in $S.$ Finally, since $C'$ is
nonnegative, we must have $C'(s)\geq \tilde{C}(s)$ for all $s \in
S.$ So we see that $C'$ contains $\tilde{C}$ and thus $D,$ so we
have already solved for $D$ by executing the moves which lead to
$C'.$

\end{proof}
As the proof suggests, numbers $n_{ij}$ can be thought of as a list
of moves, ($n_{ij}$) representing a number of moves from the $v_i$
to $v_j.$ The theorem says that if we calculate the number of
pebbles left on each vertex after the moves are all executed, by
adding a pebble for each move onto a vertex and subtracting two
pebbles for each move from it, and the resulting configuration
satisfies the demand configuration, then the original configuration
is solvable.

We may thus speak of a configuration as being \emph{solved} by a
list of moves, that is, a list of numbers $(n_{ij})$ (which we call
a \emph{solution}) meeting the condition in the theorem. We can use
this definition to expand the definition of solvability to arbitrary
(possibly negative) configurations and demand configurations saying
$C$ is \emph{solvable} for $D$ if there exists such a list of
numbers $n_{ij}$ meeting the condition of Theroem 1 for $C$ and $D.$
Clearly, this is the same as simply revising our definition of
solvability to allow negative numbers of pebbles. The theorem shows
that the two definitions concur for non-negative configurations, so
this definition generalizes the existing concept of solvability. It
will sometimes be useful to remember, however, that if a
non-negative configuration is solvable for a non-negative demand
configuration by a list of moves, it is solvable by a sequence of
pebbling moves with intermediate configurations all non-negative.

Note that Theorem 1 shows that only the difference of the two
configurations, $C-D,$ is material to whether $C$ is solvable for
$D.$ We now point out three other trivial corollaries that will be
useful in the sequel.
\begin{cor} Let $G$ be a graph with vertices $\{v_1, \ldots, v_m\}$, $C$ a configuration on $G,$ and $D$ a non-negative demand
configuration. If the list of pebbling moves $(k_{ij})$ solves $C$
for $D,$ and we have a list $(l_{ij})$ of moves with $ l_{ij}\leq
k_{ij}$ for all $1 \leq i,j \leq k$ then the configuration $C'$
obtained from $C$ by making the moves $l_{ij},$
$$C'(v_{i'}) = C(v_k) + \sum_{j'=1}^m l_{j'i'}-2\sum_{j'=1}^m l_{i'j'}$$ is
solvable by the list of moves $(n_{ij})=(l_{ij}-m_{ij})$ .
\end{cor}
That is, if we execute part of a solution, the resulting
configuration is solvable, and is solved by the list of remaining
moves.

\begin{cor}
Notation as in Theorem 1. If a configuration $C$ is solvable by
$(n_{ij})$ and $C'$ is a configuration containing $C,$ then $C'$ is
solvable by $(n_{ij})$
\end{cor}

That is, adding pebbles to a solvable configuration yields a
solvable configuration.

\begin{cor}
Notation as in Theorem 1. If $G$ has a subset of $p$ vertices,
without loss of generality $v_{m-p+1}, \ldots, v_m $ and a solvable
configuration $C$ with solution $(n_{ij})$ such that $n_{ij}=0$ if
$i \geq m-p+1$ and $j \leq m-p.$ Let $G'$ be the subgraph of $G$
induced by removing vertices $v_{m-p+1}, \ldots, v_m.$ then the
numbers $(n_{ij}),$ with $i, \, j \leq m-p$  solve the configuration
$C'$ induced by $C$ on $G'$ by restriction for the induced demand
configuration $D'.$
\end{cor}

This says that given a configuration on a graph we may remove
vertices from our graph and get a solvable configuration on the
induced subgraph if no moves from these vertices are necessary to
solve the original configuration

Next we will show that if a configuration is cover solvable, it may
be cover solved without a directed cycle of pebbling moves. This is
proven in \cite{Glenn1} using transition digraphs, but it is also a
corollary of Theorem 1. We state it here as
\begin{thm} \cite{Glenn1}
Notation as in Theorem 1. If $C$ is solvable for $D$, we can choose
the $n_{ij}$ meeting the condition in theorem 1 such that there is
no list of distinct integers $(i_1, i_2, \ldots, i_p),$ with $1 \leq
i_j \leq m$ such that $n_{i_j i_{j+1}}>0,$ for all $1\leq j < p$ and
$n_{i_p i_1} >0$.

\end{thm}
\begin{proof} Choose the list $(n_{ij})$ meeting the condition in Theorem 1
so that $\sum {n_{ij}},$ the total number of moves, is minimal. If
such integers $(i_1, i_2, \ldots, i_p)$ exist, form a new list
$(n'_{ij})$ by letting $n'_{i'j'} = n_{i'j'}-1$ if there is a $j$
such that $i' = i_j$ and $j' = i_{j+1}$ and otherwise, letting
$n'_{i'j'}=n_{i'j'}.$ Considering the new sums for each vertex
$$C(v_k)+\sum_{l=1}^m n'_{lk}-2\sum_{l=1}^m n'_{kl}.$$
we note if  $k \notin (i_1, i_2, \ldots, i_p)$ the value is
unchanged from the former sum with the $n_{ij}.$  If $k = i_j$ for
some $j$ then $ \sum n'_{lk}$ and $\sum n'_{kl}$ are both decreased
by 1 because $n'_{i_{j-1} i_j} = n_{i_{j-1} i_j} -1$ and $n'_{i_{j}
i_{j+1}} = n_{i_j i_{j+1}} -1.$ Thus the sum has been increased by 1
in this case, and so we know for all $k$ $$C(v_k)+\sum_{l=1}^m
n'_{lk}-2\sum_{l=1}^m n'_{kl}\geq C(v_k)+\sum_{l=1}^m
n'_{lk}-2\sum_{l=1}^m n'_{kl}\geq D(v_k),$$ thus the condition of
Theorem 1 is met for the $(n'_{ij}).$ But clearly
$\sum{n'_{ij}}=\sum {n_{ij}}-p$, contradicting the minimality of the
list $(n_{ij}).$
\end{proof}

This result can be applied to vertices of degree 1 in a graph in a
straightforward manner. Similar analysis is implicitly carried out
in the treatment of cover pebbling on trees in \cite{Glenn1}.

\begin{cor}
Let $G$ be a graph, $B$ a solvable configuration, $D$ a demand
configuration. Suppose $v \in G$ has degree 1, and is adjacent to
$v' \in G.$ Let $H$ be the subgraph of $G$ induced by removing $v.$
Let $D'$ be the demand function induced by restriction of $D$ to
$H.$ Construct configuration $C$ on $H$ by letting $C(w)=B(w)$ for
all $w \not=v'$ in $H$ and letting $C(v')= B(v')+ \frac{\lfloor
B(v)-D(v)\rfloor} {2}$ if $B(v)-D(v) \geq 0$ (where as usual
$\lfloor x \rfloor= \sup\{k \in \mathbb{Z} : k \leq x\}$) or $C(v')=
B(v')+ 2(B(v)-D(v))$ otherwise. Then $C$ is solvable for $D'$ on
$H.$
\end{cor}
\begin{proof}
Label the vertices of $G$ as $v_1, \, v_2, \ldots, v_m$ with $v=v_m$
and $v'=v_{m-1}$ and choose a directed-cycle-free solution
$(n_{ij})$. Since $v_m$ is adjacent only to $v_{m-1}$ we have
$n_{mj}=0$ and $n_{jm}=0$ if $j \not= m-1.$ So we must have
\begin{equation}B(v_m)+n_{(m-1) m}-2 n_{m (m-1)} \geq D(v_m).\end{equation}
Since our solution is cycle-free we must have $n_{(m-1) m}=0$ or
$n_{m (m-1)}=0$ Perform all moves $n_{(m-1) m}$ or $n_{m (m-1)}$
between $v_m$ and $v_{m-1},$ as in Corollary 2, to obtain solvable
configuration $B',$ solvable by remaining list remaining of moves.
This list has no more moves to or from $v_m$ so we apply Corollary 4
to see that the configuration $B''$ on $H$ induced by restricting
$B'$ to $H$ is solvable. Clearly we know $B''$ is equal to $B$ and
thus $C$ on every vertex of $H \setminus \{v_{m-1}\}.$

Now suppose we have $n_{(m-1) m}=0.$ By (1) we see $B(v_m)-D(v_m)
\geq 2 n_{m (m-1)}.$ Clearly $B(v_m)-D(v_m)\geq 0$ and since both
sides are integers, we have $\frac{\lfloor B(v_m)-D(v_m)\rfloor }{2}
\geq n_{m (m-1)}$ so
$$B''(v_{m-1})=B(v_{m-1})+n_{m(m-1)}\leq B(v')+\frac{\lfloor
B(v_m)-D(v_m)\rfloor }{2} =C(v').$$ Similarly, if $n_{m (m-1)}=0$ by
(1) we have $B(v_m)-D(v_m)\geq - n_{(m-1) m}$ so we have
$$B''(v_{m-1})=B(v_{m-1})-2 n_{(m-1)m}\leq B(v')+2(B(v_m)-D(v_m))
=C(v').$$ We thus know $C$ contains solvable configuration $B''$ and
so is solvable by Corollary 3.
\end{proof}

This corollary says that if a graph $G$ has solvable configuration
$C$ and a vertex $v$ of degree 1, we can just move the excess
pebbles on $v$ to its neighbor or meet the excess demand of $v$ with
pebbles from its neighbor and then get rid of $v,$ and the resulting
configuration on the induced subgraph will be solvable. In
particular, this result allows us to easily tell if a configuration
on a tree is solvable.

There is one more theorem we need in our toolbox. Given graph $G,$
configuration $C$ and demand function $D$ define function
$\Gamma_C(v)$ by
$$\Gamma_C(v)= \sum_{v' \in G} (C(v')-D(v'))2^{-d(v',v)}$$ where $d$
represents graph theoretic distance.

\begin{thm}
Let $G$ be a graph, $C$ a configuration, $D$ a demand configuration.
If $\ \Gamma_C(v)< 0$ for some $v \in G,$ then $C$ is not solvable.
\end{thm}
\begin{proof}
Clearly $\Gamma_C(v)$ can not be increased by a pebbling move, since
a pebbling move always removes two pebbles from a vertex and adds
one pebble to a vertex which is at most 1 edge closer to $v.$ A
solved configuration has $C(v') \geq D(v')$ for all $v'\in G$ thus
$\Gamma_C(v')\geq 0$ for all $v'.$ Thus, we cannot solve a
configuration $C$ such that $\Gamma_C(v)<0$ for some $v \in G$
through a sequence of pebbling moves.\end{proof}

\section{The $N\!P$-Completeness of the Cover Solvability Problem}

The material in this section appears in \cite{monthlypaper}. I
repeat it here to give a more rigorous proof, using the tools
developed above.

\begin{thm}
 The cover solvability decision problem which accepts pairs $\{ G, C, D\}$ if and only if $G$ is a graph and $C$
 is a non-negative configuration on $G$ which is solvable for configuration $D$ is in $N\! P.$
\end{thm}

\begin{proof} Given the list of integers $n_{ij},$ the solvability
of $C$ for $D$ may be quickly checked by checking that the $|G|$
inequalities in Theorem 1 hold.
 \end{proof}

To show the problem is $N\!P$-hard, we make use of the following
problem:

\begin{defn}
Let the {\it exact cover by 4-sets problem} be the decision problem
which takes as input a set $S$ with $4n$ elements and a class $A$ of
at least $n$ four-element subsets of $S,$ accepting such a pair if
there exists an $A' \subseteq A$ such that $A'$ is a class of
disjoint subsets which make a partition of $S,$ that is they are $n$
subsets containing every element of $S.$
\end{defn}

This problem is $N\!P$-complete \cite{Karp}. Indeed, the
corresponding problem of exact cover by 3-sets is also
$N\!P$-complete, but for our purposes, the 4-set problem is more
convenient. We now show the main result of this section:

 \begin{thm}
 The cover solvability decision problem which accepts pairs $\{ G, C\}$ if and only if $G$ is a graph and $C$
 is a configuration which is solvable for the unit configuration on $V(G)$ is  $N\!
 P$-hard and thus $N\!P$-complete.
 \end{thm}

\begin{proof}We will show that instances the exact cover by 4-sets may be
translated to equivalent instances of cover solvability for the unit
configuration in polynomial time. Given an instance this problem,
that is, a set $S=\{s_1,s_2,\ldots, s_{4n}\}$ and a class $A= \{a_1,
a_2, \ldots, a_m\}$ of four-element subsets of $S,$ construct a
graph $G'$ in the following manner: create a set of vertices
$T=\{t_1, t_2, \ldots, t_{4n}\}$ corresponding to the elements of
$S,$ and a set of vertices $B = \{b_1, b_2, \ldots, b_m\}$
corresponding to the members of $A.$ Create edges between $B$ and
$T$ in the natural way, including $\{b_i, t_j\}$ if $s_j \in a_i.$
Additionally, create a vertex $v$ and a path of length $m-n$ which
has one terminal vertex $v$ and the other called $w.$ Finally,
create vertex classes $ B' = \{{b_1}',{b_2}', \ldots, {b_m}' \} $
and $B''= \{{b_1}'',{b_2}'', \ldots, {b_m}'' \},$ creating edges
$\{b_i, {b_i}'\},$ $\{{b_i}',{b_i}''\}$ and $\{{b_i}'', v\}$  for
all $i.$

Now we create a configuration $C'$ on $G'.$ Let $C'(t)=0$ for all $t
\in T$ and let $C'(b) =9$ for all $b \in B.$ Let $C'(v) =
2^{m-n}-(m-n)+1,$ $C'(w)=0$ and $C'(u)=1$ for all $u$ between $v$
and $w$ on the path connecting them. Let $C'(u) = 1 $ for all $u \in
B' \cup B''.$ (Figure \ref{dia}.)

 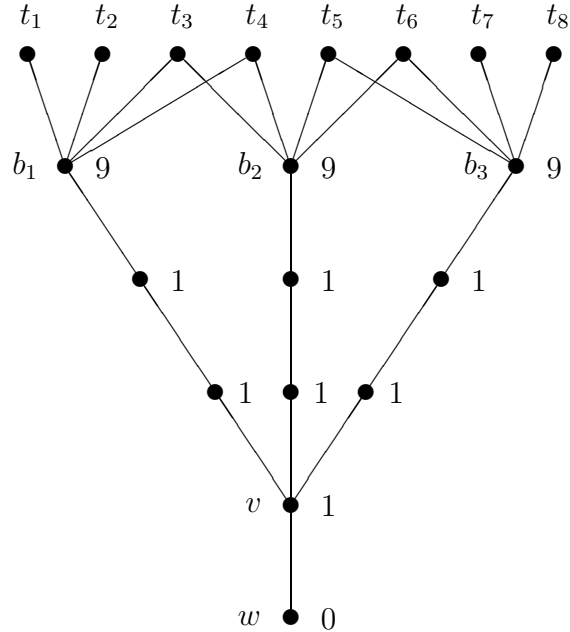
\begin{figure}[htb]
\unitlength 1mm
\begin{center}
\begin{picture}(80,90)
%root of the tree
\put(40,20){\circle*{2}} \put(40,35){\circle*{2}}
\put(30,35){\circle*{2}} \put(50,35){\circle*{2}}
\put(40,20){\line(-2,3){10}} % \put(x,y){\line(xslope,yslope){length}
\put(40,20){\line(0,1){15}} \put(40,20){\line(2,3){10}} %level 1
\put(40,50){\circle*{2}} %the * makes the circle darkened
\put(60,50){\circle*{2}} \put(20,50){\circle*{2}}
\put(40,50){\line(0,1){15}} % \put(x,y){\line(xslope,yslope){length}
\put(60,50){\line(2,3){10}}
\put(20,50){\line(-2,3){10}} % \put(x,y){\line(xslope,yslope){length}
\put(40,35){\line(0,1){15}} % \put(x,y){\line(xslope,yslope){length}
\put(50,35){\line(2,3){10}}
\put(30,35){\line(-2,3){10}} % \put(x,y){\line(xslope,yslope){length}
\put(40,65){\circle*{2}} %the * makes the circle darkened
\put(70,65){\circle*{2}} \put(10,65){\circle*{2}}
\multiput(5,80)(10,0){8}{\circle*{2}} \put(4, 84){$t_1$}
\put(14,84){$t_2$}\put(24, 84){$t_3$} \put(34,84){$t_4$} \put(44,
84){$t_5$} \put(54,84){$t_6$}\put(64, 84){$t_7$} \put(74,84){$t_8$}
%the * makes the circle darkened
%level 2
\put(40,65){\line(-1,1){15}} \put(40,65){\line(-1,3){5}}
\put(40,65){\line(1,3){5}} \put(40,65){\line(1,1){15}}
\put(10,65){\line(-1,3){5}} \put(10,65){\line(1,3){5}}
\put(10,65){\line(1,1){15}} \put(10,65){\line(5,3){24}}
\put(70,65){\line(1,3){5}} \put(70,65){\line(-1,3){5}}
\put(70,65){\line(-1,1){15}} \put(70,65){\line(-5,3){24}}
%labels
\put(14,63.5){9} \put(44,63.5){9} \put(74,63.5){9} \put(24,48.5){1}
\put(44, 48.5){1} \put(64, 48.5){1} \put(33,33.5){1} \put(43,
33.5){1} \put(53, 33.5){1} \put(3,64){$b_1$}
\put(33,64){$b_2$}\put(63,64){$b_3$} \put(34, 19){$v$}
\put(44,18.5){1}

\put(40,5){\line(0,1){15}} \put(40,5){\circle*{2}} \put(33, 4){$w$}
\put(44,3.5){0}
\end{picture}
\end{center}
\caption{\label{dia} A cover solvability problem that corresponds to
the exact cover by four $4$-sets problem, $a_1 =
\{s_1,s_2,s_3,s_4\}$, $a_2 = \{s_3,s_4,s_5,s_6\}$, $a_3 = \{
s_5,s_6,s_7,s_8 \}$.} \label{ex2}
\end{figure}

Since we have made only $5|A|+1$ vertices, the construction may be
done in polynomial time. To finish the proof, we now must that $C'$
is solvable if and only if $A$ contains an exact cover of $S.$

First suppose that $A$ contains a exact cover $A' = \{ a_{i_1},
a_{i_2}, \ldots, a_{i_n} \} $ of $S.$ Then for each vertex in $B$
which is a $b_{i_j}$ for some $ 1 \leq j \leq n,$ we use 8 of the
pebbles on this vertex to put one pebble on each of the four
vertices of $T$ to which it is adjacent. Because of the fact that
$A'$ is a exact cover and the way we constructed $G',$ we now have
one pebble on every vertex of $T.$ Furthermore, we have $m-n$
vertices in $B$ that still have 9 pebbles each on them. Because $v$
is at distance $3$ from each of these vertices, we can use 8 pebbles
from each of these vertices to move one pebble each onto $v,$
leaving $2^{m-n} +1$ pebbles on $v,$ enough to move one pebble onto
$w$ while leaving one pebble on $v.$ This leaves this exactly one
pebble on every vertex of $G',$ so we have solved $C'$ for the unit
configuration. Therefore, $C'$ is solvable whenever $A$ contains a
exact cover of $S.$

To show the converse, suppose that $A$ does not contain a exact
cover of $S.$ We must show $C'$ is not solvable for the unit
configuration. Assume the opposite, $C'$ is solvable. The sequence
of pebbling moves which solves $C'$ must contain (at least) one move
to $t$ for every $t \in T.$ Of the $4n$ moves that are necessary
(one to each vertex in $T,$) all must originate from $B,$ and no
more than 4 can originate from any one vertex of $B,$ since each
vertex in $B$ is adjacent to only four vertices in $T.$ Since $A$
does not contain a exact cover of $S,$ it cannot be the case that
these moves originate from exactly $n$ vertices in $B.$

We make these $4n$ moves immediately from $C',$ using Corollary 2 to
see that the resulting configuration must be solvable. There are now
$8(m-n)+m$ pebbles left on $B,$ total. Suppose another move to $B$
is necessary to solve the graph. Then we can make one such move
immediately. But then we apply Theorem 7 to show that the resulting
configuration is unsolvable, since
$$\Gamma_U(w)=2^{-(m-n)}(2^{m-n}-(m-n))+2^{-(m-n+3)}(8(m-n)-2)+2^{-(m-n+4)}<0$$ can be easily seen. We
now know that none of the remaining moves are onto vertices in $T,$
and thus there can be no more moves from $T$ because such a move
would leave a vertex in $T$ permanently uncovered. Our configuration
is solvable without any further moves to or from any vertex in $T.$
We can thus eliminate these vertices $T$ from our graph using
Corollary 4, leaving us with a tree.

Corollary 6 now allows us to collapse the paths of length three
which terminate in the vertices of $B,$ and we easily see that
collapsing these paths one vertex at a time adds one pebble to $v$
if there are still $9$ pebbles on the end vertex in $B,$ and adds
zero pebbles to $v$ otherwise. But since it is not the case that the
moves originate from exactly $n$ vertices in $B,$ there are less
than $m-n$ total stacks of $9$ pebbles on $B$ left, and so we add
less than $m-n$ pebbles to $v$ through this process. This leaves us
with less than $2^{m-n}+1$ pebbles on $v,$ and therefore we do not
have enough pebbles to move a pebble onto $w$ while leaving a pebble
on $v.$ Therefore, the configuration on this reduced graph is
unsolvable, so the original configuration was unsolvable,
contradiction.
\end{proof}

\section{The $N\!P$-Completeness of the Canonical Pebbling Solvability Problem}

A configuration $C$ on a graph $G$ is canonical pebbling solvable if
and only if the $v$-reachability configuration $R_v$ is solvable for
all $v \in G$ Thus, the canonical pebbling solvability question is
equivalent to $|G|$ cover pebbling questions, and is therefore
trivially in $N\!P.$ We now show it is in fact also $N\!P$-hard and
thus is $N\!P$-complete.

\begin{thm}
 The traditional pebbling solvability decision problem, which accepts pairs $\{ G, C\}$ if and only if $G$ is a graph and $C$
 is a configuration which is solvable for each configuration which has value 1 on a single vertex in $V(G)$ and value 0 on every other vertex
 is $N\!P$-hard and thus $N\!P$-complete.
\end{thm}

\begin{proof}We demonstrate that it is possible to ``translate''
a cover pebbling solvability question into a pebbling solvability
question in polynomial time, thus showing the pebbling solvability
question is $N\!P$-hard since cover solvability is $N\!P$-hard by
Theorem 10. For simplicity, we translate instances of the restricted
class of problems shown $N\!P$ hard by Theorem 10, that is, problems
of the solvability of non-negative configurations for the unit
configuration.

Consider a graph $G$ with $|G|=n$ and a non-negative configuration
$C$ on $G$ to be solved for the unit configuration. If $|C| \geq
2^n-1$ then $C$ is solvable because we know in \cite{Glenn1} that
$2^n-1$ is the largest possible cover pebbling number possible for
the unit configuration on a graph with $n$ vertices. So if $|C| \geq
2^n$ we translate our solvability problem to some trivially solvable
canonical pebbling solvability problem, such as the solvability
problem of the empty graph. We may therefore assume that $|C| <
2^n.$ Further, if $|G|=1,$ it is obvious that a configuration is
canonical pebbling solvable if and only if it is cover solvable for
the unit configuration, so we leave the graph and configuration
unchanged by our translation. So we may assume $n \geq 2.$

Label the vertices of $G$ $\{v_1,v_2, \ldots, v_n\}.$ Construct a
graph $G'$ and a configuration $C'$ in the following manner: begin
with an isomorphic copy of $G$  with vertices $\{{v_1}',
{v_2}',\ldots, {v_n}'\}$ and with edges added in the obvious way,
that is with $\{{v_i}',{v_j}'\} \in E(G')$ if and only if
$\{v_i,v_j\} \in E(G)$. Let $H$ be the subgraph induced by
$\{{v_1}', {v_2}',\ldots, {v_n}'\},$ which is isomorphic to $G.$ Now
we add a path of length $n$ to each vertex ${v_i}'$ and attach the
opposite end of each path to a vertex $w_0$. To put this precisely,
add vertices $u_{ij}$ for all $1\leq i, \ j \leq n,$ and a vertex
$w_0$ and edges $\{{v_i}', u_{i1}\},$ $\{u_{ij},u_{i(j+1)}\},$ and
$\{u_{in},w_0\}$ for all $1 \leq i \leq n$ for all $1 \leq j \leq
n-1.$ Finally, add a path of length $n$ attached at $w_0$ and call
the terminal vertex $w_n.$ That is, add vertices $\{w_1,w_2, \ldots,
w_n \}$ and edges $\{w_{i-1}, w_i\}$ for each $1 \leq i \leq n.$

Now, we define a configuration $C'$ on the graph $G'$ we
constructed: let $C'({v_i}')=C(v_i)+1$ for all $1 \leq i \leq n.$
Let $C'(u_{ij}) = 1$ for all $1 \leq i, \ j \leq n.$ Let $C'(w_0) =
2^n - n$ and let $C'(w_i) = 0$ for all $1 \leq i \leq n.$ (Figure 2.)

 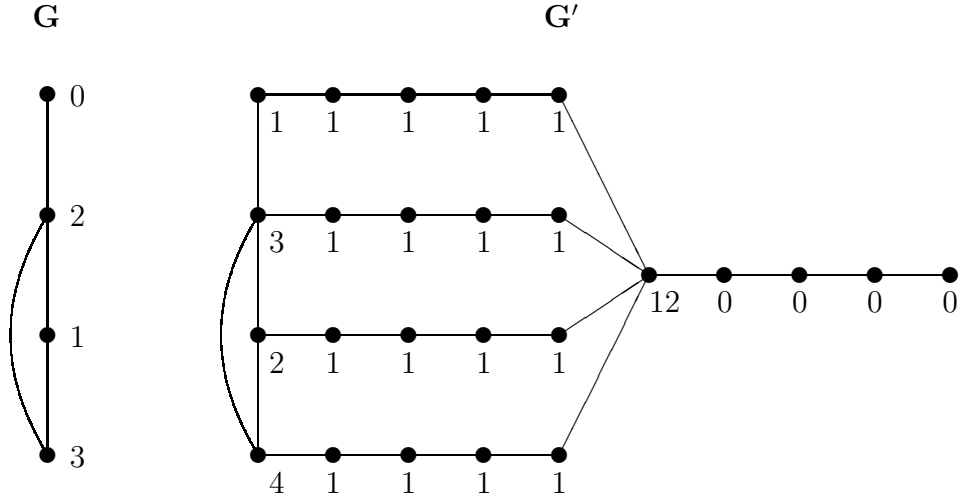
\begin{figure}[htb]

\unitlength 1mm
\begin{center}
\begin{picture}(140,76)

\put(78,66){$\mathbf{G'}$} \put(10,66){$\mathbf{G}$}

%graph G
\multiput(12,57)(0,-16){4}{\circle*2}\multiput(12,57)(0,-16){3}{\line(0,-1){16}}
\qbezier(12,41)(2,25)(12,9)

%four parallel paths:
\multiput(40,57)(10,0){5}{\circle*{2}} \multiput(40,
57)(10,0){4}{\line(1,0){10}}

\multiput(40,41)(10,0){5}{\circle*{2}} \multiput(40,
41)(10,0){4}{\line(1,0){10}}

\multiput(40,25)(10,0){5}{\circle*{2}} \multiput(40,
25)(10,0){4}{\line(1,0){10}}

\multiput(40,9)(10,0){5}{\circle*{2}} \multiput(40,
9)(10,0){4}{\line(1,0){10}}

%lines on the left
\multiput(40,57)(0,-16){3}{\line(0,-1){16}} \qbezier(40, 41) (30,
25) (40, 9)

%the lines to w_0

\put (80,57){\line(1,-2){12}} \put (80,41){\line(3,-2){12}}
\put(80,25){\line(3,2){12}}  \put (80,9){\line(1,2){12}}

%path on right

\multiput(92,33)(10,0){5}{\circle*{2}} \multiput(92,
33)(10,0){4}{\line(1,0){10}}

%labels of G
\put(15,55.5){0} \put(15, 39.5){2} \put(15, 23.5){1} \put(15,
7.5){3}

%labels of G'
\put(41.5,52){1} \put(41.5, 36){3} \put(41.5, 20){2} \put(41.5,
4){4}

\multiput(49, 52)(10,0){4}{1} \multiput(49, 36)(10,0){4}{1}
\multiput(49, 20)(10,0){4}{1} \multiput(49, 4)(10,0){4}{1}

\put(92,28){12} \put(101,28){0} \put(111,28){0} \put(121, 28){0}
\put(131, 28){0}

\end{picture}
\end{center}
\caption{\label{dia2} Translating a cover solvability problem on $G$
to a canonical pebbling solvability problem on the graph $G'$}
\end{figure}

We have created $|G|^2 + 2|G| + 1$ vertices, so this construction
may be done in polynomial time. We now show that the configuration
$C'$ is canonical pebbling solvable on $G'$ if and only if $C$ is
solvable for the unit configuration on $G.$

We know that if the $w_n$-reachability configuration on $R_{w_n}$ is
not solvable, then $w_n$ is not reachable from $C'$ and so $C'$ is
not canonical pebbling solvable. Also if $R_{w_n}$ is solvable from
$C'$, we know from Theorem 1 it may be solved by a sequence of moves
for which the intermediate configurations are all non-negative.
Clearly, this requires moving pebbles onto each vertex $w_1, \ldots,
w_{n-1}$ in some intermediate configuration, thus these vertices are
reachable from $C'.$ Since every other vertex of $G'$ is trivially
reachable from $C',$ we know $C'$ is canonical pebbling solvable.
Thus the canonical pebbling solvability of $C'$ is equivalent to the
solvability of $R_{w_n}.$ Therefore, we henceforth consider the
solvability of $R_{w_n}.$

If the configuration $C$ is solvable for the unit configuration on
$G,$ we may perform the same moves that solve $C$ on the isomorphic
copy of $H$ in $G'$ and because each vertex in $H$ begins with one
more pebble on each vertex, the result will be a configuration with
one more pebble on each vertex of $H$ than is on the corresponding
vertex of $G$ after the configuration $C$ has been solved. In
particular, there will be at least 2 pebbles on each vertex
${v_i}'$. Thus, we can move a pebble from each ${v_i}'$ onto
$u_{i1},$, which leaves 2 pebbles on each vertex $u_{i1},$ meaning
we may move a pebble from each $u_{i1}$ onto each $u_{i2},$ and so
forth, until we have 2 pebbles on each $u_{in},$ from each of which
we may move a pebble onto $w_0,$ adding $n$ pebbles to $w_0$ for a
total of $2^n.$ We can use these $2^n$ pebbles to add one pebble to
$w_n,$ which is at distance $n$ from $w_0.$ So $R_{w_n}$ is solvable
whenever $C$ is solvable for the unit configuration on $G.$

Now we turn to the converse. Assume the converse does not hold, that
$C$ is not solvable for the unit configuration on $G,$ but $C'$ is
solvable for $R_{w_n}$ on $G'.$ Then, using Theorem 5, $C'$ has a
directed-cycle-free solution for $R_{w_n},$ so choose such a
solution. Suppose there is a move from $u_{i1}$ to ${v_i}'$ for some
$i$. Since the move from $u_{i1}$ to ${v_i}'$ leaves $-1$ pebbles on
$u_{i1},$ we must have a move onto $u_{i1}.$ But our solution is
directed-cycle-free, so we do not have a move from ${v_i}'$ to
$u_{i1},$ so the only possibility is a move from $u_{i2}$ to
$u_{i1}.$ By the same reasoning, we must have a move onto $u_{i2}$
which cannot originate from $u_{i1}$, and thus must originate from
$u_{i3}.$ Repeating the argument, we have a move from $u_{i1}$ to
$v_i,$ a move from $u_{ij}$ to $u_{i(j-1)}$ for $2 \leq j \leq n,$
and a move from $w_0$ to $u_{in}.$ By Corollary 2, the configuration
$C''$ obtained by executing all of these moves is solvable.

We now use Theorem 7 to derive a contradiction, by calculating
$\Gamma_{C''}(w_n).$ We have added one pebble to $H,$ so there are
now at most $2^n$ pebbles on $H.$ Since these vertices have distance
$2n+1$ from $w_n$ this gives a contribution to $\Gamma_{C''}(w_n)$
of at most $2^n 2^{-2n+-1}=2^{-n-1}.$ For the $u_{ij}$ we see $n-1$
paths are still covered with pebbles, so we calculate a contribution
to $\Gamma_{C''}(w_n)$ of $$(n-1)\sum_{j=n+1}^{2n} 2^{-j} =
(n-1)(2^{-n} - 2^{-2n})= n 2^{-n}-2^{-n}-(n-1)2^{-2n}.$$ Finally,
the contribution of $w_0$ is $(2^n - n) 2^{-n} = 1 - n 2^{-n}$ and
the contribution of $w_n$ is $-1.$ Totalling up, we see
$$\Gamma_{C''}(w_n) \leq -1 +(1 - n 2^{-n})+
(n2^{-n}-2^{-n}-(n-1)2^{-2n})+ 2^{-n-1} = -2^{-n+1}+2^{-2n}$$ and
since $n \geq 2$ we have $\Gamma_{C''}(w_n)<0$ so the $C''$ is not
solvable for $R_{w_n},$ contradiction.

Therefore, there can be no move from $u_{i1}$ to ${v_i}',$ for any
$i.$ Now execute every move on our solution from some ${v_i}'$ to
$u_{i1},$ and call the resulting configuration $C''.$ By Corollary
2, the resulting configuration is solvable by the remaining moves.
But there are no moves left from $H$ to $G'\setminus H$ and no moves
from $G' \setminus H$ to $H$, so by Corollary 4, the configuration
induced by $C''$ on $H$ and the configuration induced by $C''$ on
$G' \setminus H$ must both be solvable. Consider first the
configuration induced by $C''$ on $G' \setminus H.$ $G' \setminus H$
is a tree, and using Corollary 6, we may collapse the paths formed
by the $u_{ij}$ down to $w_0.$ We know that $C''$ agrees with $C'$
on $G' \setminus H$ except possible on the vertices $u_{i1}$ for
each $i,$ and we easily see that when we collapse the path formed by
$\{u_{i1}, u_{i2}, \ldots, u_{in}\}$ to $w_0$ we add no pebbles to
$W_0$ if $C''(u_{i1}) \leq 1$ and one pebble if $2 \leq
C''(u_{i1})\leq 2^{n-1} + 1.$

First suppose $C''(u_i1) \geq 2^{n-1}+2$ for some $i.$ This is
possible only if we have made at least $2^{n-1}+1$ moves from $v_i$
to $u_{i1}.$ Now consider the restriction of $C''$ to $H$. We had at
most $2^n-1$ pebbles on $H$ originally, and the moves from $v_i$ to
$u_{i1}$ decrease this number by at least $2^n+2.$ We thus have at
most $-3$ pebbles left on $H,$ and since pebbling moves can only
decrease this number, we can never get back to 0, which is required
to solve for $R_{w_n}$ restricted to $H.$ So the restriction of
$C''$ to $H$ is not solvable for this restricted demand function,
which is impossible by Corollary 6.

Now suppose $C''(u_{i1})\leq 2^{n-1} + 1$ for all $i.$ Then
collapsing the paths formed by $\{u_{i1}, u_{i2}, \ldots, u_{in}\}$
to $w_0$ we add one pebble to $w_0$ for each $i$ such that
$C''(u_{i1})\geq 2.$ If this is not every $i,$ we are left with less
than $2^n$ pebbles on $w_0,$ which clearly makes $w_n$ unreachable,
so we may assume $C''(u_{i1})\geq 2$ for each $i.$ This means we
have made at least one move from ${v_i}'$ to $u_{i1}$ for each $i,$
decreasing the number of pebbles on each ${v_i}$ by at least 2. We
know by Corollary 6 that the resulting configuration $C''$
restricted to $H$ is solvable for $R_{w_n}$ restricted to $H,$ that
is, the configuration which is identically 0 on $H.$ Since we know
$C''({v_i}') \leq C'({v_i}')-2 = C(v_i)-1$ for all $i$ we know by
Corollary 3 that the configuration $C' - 2$ found by starting with
$C'$ and removing two pebbles from each vertex of $H$ is solvable
for the zero configuration on $H.$ By the isomorphism of $G$ and $H$
and the fact that $C'({v_i})-2= C(v_i)-1$ for all $i,$ this is
equivalent to saying $C(v_i)-U$ is solvable for the zero
configuration. But by Theorem 1, we know that only the difference
between the demand function and the initial configuration is
material to cover solvability. Thus, this is equivalent to saying
$C$ is solvable for the unit configuration on $G,$ which we assumed
is false, contradiction.
\end{proof}

\begin{cor} The decision problem which accepts pairs $\{ G, C, v\}$ if and only
if $G$ is a graph, $v$ a vertex of $G$ and $C$ is a non-negative
configuration on $G$ which is solvable for configuration $R_v$ is
$N\! P$-complete.
\end{cor}
\begin{proof} We showed above that the canonical pebbling solvability question we constructed
is equivalent to problem of solvability for the configuration
$R_{w_n}.$ Thus, the proof also shows that the cover pebbling
solvability problem restricted to reachability configurations is
still $N\!P$-hard, and thus $N\!P$-complete.
\end{proof}

\section{The $N\!P$-Hardness of Determining Pebbling Numbers for Reachability
Configurations}

In \cite{Vuong} a simple formula for the cover pebbling number of a
strictly positive configuration on a graph is given. In this
section, however, the problem of determining whether the pebbling
number of a general non-negative configuration is greater than a
certain quantity is shown to be $N\!P$-hard.

\begin{thm}
The decision problem which accepts triples $\{G, v, m\}$ if and only
if $\gamma_G(R_v) > m$ is $N\!P$-hard.
\end{thm}
\begin{proof}
Again, we use the exact cover by four-sets problem. Given a set
$S=\{s_1,s_2,\ldots, s_{4n}\}$ and a class $A= \{a_1, a_2, \ldots,
a_m\}$ of four-element subsets of $S,$ construct a graph $G'$ as
follows: as before create a set a set of vertices $T=\{t_1, t_2,
\ldots, t_{4n}\}$ corresponding to the elements of $S,$ and a set of
vertices $B = \{b_1, b_2, \ldots, b_m\}$ corresponding to the
members of $A,$ and include edge $\{b_i, t_j\}$ whenever $s_j \in
a_i.$ Create a vertex $v$ and add edges $\{t_i,v\}$ for all $1 \leq
i \leq 4n.$ Finally, attach a path of length 3 to each vertex in
$B,$ that is, create vertices ${b_j}', {b_j}'',$ and ${b_j}'''$ and
edges $\{b_j,{b_j}'\},$ $\{{b_j}',{b_j}''\},$ and
$\{{b_j}'',{b_j}'''\},$ for each $1 \leq j \leq m.$ (Figure 3.)

\begin{figure}[htb]
\unitlength 1mm
\begin{center}
\begin{picture}(80,85)

%top three levels of dots
\put(40,80){\circle*{3}} \put(35,79){$v$}

 \multiput(5,65)(10,0){8}{\circle*{2}}

\put(40,50){\circle*{2}} \put(70,50){\circle*{2}}
\put(10,50){\circle*{2}}

%the slanty lines at the top

\qbezier[140](5,65)(23.5,72.5)(40,80)
\qbezier[140](75,65)(57.5,72.5)(40,80) \put(15, 65){\line(5,3){25}}
\put(25,65){\line(1,1){15}} \put(35,65){\line(1,3){5}}
\put(45,65){\line(-1,3){5}} \put(55,65){\line(-1,1){15}}
\put(65,65){\line(-5,3){25}}

%the slanty lines second from the top
\put(40,50){\line(-1,1){15}} \put(40,50){\line(-1,3){5}}
\put(40,50){\line(1,3){5}} \put(40,50){\line(1,1){15}}
\put(10,50){\line(-1,3){5}} \put(10,50){\line(1,3){5}}
\put(10,50){\line(1,1){15}} \put(10,50){\line(5,3){24}}
\put(70,50){\line(1,3){5}} \put(70,50){\line(-1,3){5}}
\put(70,50){\line(-1,1){15}} \put(70,50){\line(-5,3){24}}

%the paths at the bottom

\multiput(10,50)(0,-15){4}{\circle*{2}}
\multiput(40,50)(0,-15){4}{\circle*{2}}
\multiput(70,50)(0,-15){4}{\circle*{2}}

\multiput(10,50)(0,-15){3}{\line(0,-1){15}}
\multiput(40,50)(0,-15){3}{\line(0,-1){15}}
\multiput(70,50)(0,-15){3}{\line(0,-1){15}}

%labels
\put(3.5,49){$b_1$} \put(33.5,49){$b_2$}\put(63.5,49){$b_3$}
\put(3,34){${b_1}'$} \put(33,34){${b_2}'$}\put(63,34){${b_3}'$}
\put(2.5,19){${b_1}''$}
\put(32.5,19){${b_2}''$}\put(62.5,19){${b_3}''$}
\put(2,4){${b_1}'''$} \put(32,4){${b_2}'''$}\put(62,4){${b_3}'''$}

\put(0, 64){$t_1$} \put(10,64){$t_2$}\put(20, 64){$t_3$}
\put(30,64){$t_4$} \put(40, 64){$t_5$} \put(50,64){$t_6$}\put(60,
64){$t_7$} \put(70,64){$t_8$}

\end{picture}
\end{center}
\caption{\label{dia3} A cover pebbling number problem that
corresponds to the exact cover by four $4$-sets problem, $a_1 =
\{s_1,s_2,s_3,s_4\}$, $a_2 = \{s_3,s_4,s_5,s_6\}$, $a_3 = \{
s_5,s_6,s_7,s_8 \}$.} \label{ex2}
\end{figure}
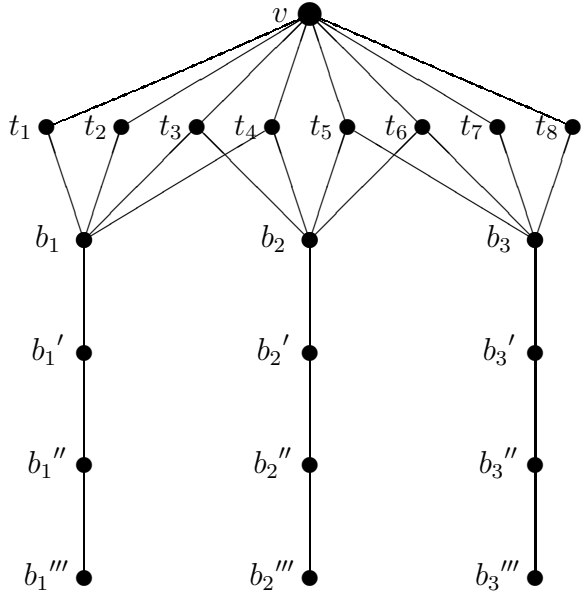

Now, we claim that $\gamma_{G'}(R_v) > 31 n + 15 (m-n)= 15 m + 16 n$
if and only if $A$ contains an exact cover of $S.$ First, suppose
such an exact cover $\{a_{i_1}, \ldots, a_{i_n}\}$ exists. Define
configuration $C$ on $G$ by putting 31 pebbles on each
$\{{b_{i_1}}''', \ldots, {b_{i_n}}'''\}$ and 15 on each ${b_i}'''$
such that $i \not=i_j$ for all $j.$ We easily calculate
$|C|=31n+15(m-n).$ We claim this configuration is not solvable for
$R_v$, and so $\gamma_G(R_v) > 31 n + 15 (m-n).$

Applying Corollary 6, we can collapse the paths to the $b_i,$
leaving 2 pebbles on $b_{i_j}$ for all $j$ and 1 pebble on all $b_i$
for which $i \not= j$ for all $j.$ But since $\{a_{i_1}, \ldots,
a_{i_n}\}$ is an exact cover, the $b_{i_j}$ are not mutually
adjacent to any vertex, so while we may move the stacks of 2 onto
vertices of $T,$ we can only produce stacks of at most one pebble on
any vertex, so no further moves are possible, and we cannot solve
for $R_v.$

Now suppose no exact cover exists. We wish to show an arbitrary
configuration  $C'$ of size $31n+15(m-n)$ on $G'$ is solvable for
$R_v$ and thus $\gamma_{G'}(R_v) \leq 31n+15(m-n).$ Clearly, any
configuration which has a pebble on $v$ is trivially solvable for
$R_v,$ so we may assume our configuration has no such pebble.
Similarly, if our configuration has more than two pebbles on any
vertex in $T,$ it is solvable for $R_v$ by using these two pebbles
to put a pebble on $v.$ So we may assume there is at most one pebble
on any vertex of $T.$

We consider two cases. First suppose the configuration has no
pebbles on $T$. Then all $16n+15m$ pebbles are distributed on the
$m$ paths $(b_{i}, {b_i}', {b_i}'', {b_i}''').$ Note that we know
from \cite{Glenn2} that the canonical pebbling number of a path of
length 4 is 8. Thus, we can use any collection eight pebbles on one
of these paths to move a pebble to $b_{i}.$ If there are 32 pebbles
or more on any of the paths $(b_{i}, {b_i}', {b_i}'', {b_i}'''),$ we
may use them to move 4 pebbles to $b_i,$ and from there two onto a
vertex in $T,$ and then one onto $v,$ so we may assume we have no
more than 31 on any of the paths. Thus, we must have 16 pebbles on
at least $n$ of the paths, because any configuration which avoids
this has at most $31(n-1)+15(m-(n-1))= 16n+15m-16$ pebbles
distributed on the paths, fewer than the $16n+15m$ in the
configuration we are considering. We can use these $n$ collections
of 16 pebbles to move two pebbles each onto $n$ of the vertices in
$B.$ Since there is no exact cover of $S$ in $A,$ we know that of
this set of $n$ of vertices in $B,$ two are mutually adjacent to a
vertex in $T$. We can move one pebble each from these two vertices
onto the vertex in $T$ which they are both adjacent to, and then use
those two pebbles to put one pebble on $v.$ Thus, the original
configuration must be solvable for $R_v$.

Now suppose there $k>0$ pebbles on $T.$ If there are more than 15
pebbles on a path $(b_{i}, {b_i}', {b_i}'', {b_i}''')$ with $b_i$
adjacent to a vertex in $T$ which has a pebble on it, we can move
two pebbles onto $b_i$ using these 16 or more pebbles, use these to
add a pebble to this vertex of $T,$ and then use these two pebbles
to put a pebble on $v,$ showing our configuration is solvable for
$R_v.$ So we assume none of these paths has 16 or more pebbles.
Also, as before, we may assume none of the $m$ paths below $B$ has
more than 31 pebbles total. Now note that any set of $n -
\left\lfloor \frac{k-1}{4}\right\rfloor$ of the vertexes in $B$
which are not adjacent to any of the $k$ covered vertices in $T$
contains two vertices which are mutually adjacent to some vertex in
$T$ by the pigeonhole principle (there are $4n-k$ vertices in $T$
which are not covered, and every vertex in $B$ which is not adjacent
to a covered vertex is adjacent to four non-covered vertices.)

We now claim that we have at least 16 pebbles on $n - \left\lfloor
\frac{k-1}{4}\right\rfloor$ of the paths below $B$ which are not
rooted at to any of the vertices in $B$ which are adjacent to
covered vertices in $T.$ If our configuration avoided having $n -
\left\lfloor \frac{k-1}{4}\right\rfloor$ such paths with 16 pebbles,
we calculate the maximum number of pebbles on the graph as
$31\left(n - \left\lfloor \frac{k-1}{4}\right\rfloor -
1\right)+15\left(m-\left(n - \left\lfloor
\frac{k-1}{4}\right\rfloor-1\right)\right)+k=15m+16n+k-16\left\lfloor
\frac{k-1}{4}\right\rfloor-16<15m+16n,$ less than the number of
pebbles in our configuration. So we have 16 pebbles on $n -
\left\lfloor \frac{k-1}{4}\right\rfloor$ of the paths below $B$
which are not adjacent to a covered vertex in $T,$ and using the
above, we see we can move two pebbles to each of two vertices in $B$
which are mutually adjacent to a vertex in $T.$ Move one pebble from
each of these vertices to this vertex of $T,$ and then use these two
pebbles to place one on $v.$ We have now shown $R_v$ is solvable in
all cases if our configuration is of size $31n+15(m-n).$ Thus
$\gamma_{G'}(R_v) \leq 31n+15(m-n)$ if $A$ contains no exact cover
of $S.$ We now know our cover pebbling number problem is equivalent
to the exact cover problem. Thus the problem of determining in
general the cover pebbling number of a reachability configuration is
$N\!P$-hard.
\end{proof}

\paragraph*{Open Questions}
\begin{trivlist}
\item \textbf{1}. Does there exist highly
symmetric class of graphs such that solvability questions on this
class of graphs are still $N\!P$-complete? For instance, are cover
solvability questions on the hypercube graph $N\!P$ complete?

\item \textbf{2}. (Suggested in \cite{Milans}.) A graph $G$ for which $\pi(G) = |G|$ has been
called \mbox{\emph{Class 0}}. These graphs have been the focus of
much research, and an interesting sufficient condition for a graph
to be \mbox{Class 0} is given in \cite{Glenn4}. What is the
complexity of the problem which asks whether a graph is \mbox{Class
0}?

\item \textbf{3}. (Also suggested in \cite{Milans}.) What is the
complexity of the problems discussed in this paper when the graphs
are restricted to be planar? Outerplaner?

\end{trivlist}

\paragraph*{Acknowledgments}

Part of this work was done while the author was receiving support
from NSF grant DMS-0139286 while a participant at the East Tennessee
State University REU under the direction of Anant Godbole. I
acknowledge him for his guidance and encouragement. Also, I thank
Carl Yerger for his help preparing the original version of Section 3
for publication in \cite{monthlypaper}, including providing comments
and creating Diagram 1. Finally, I acknowledge Kevin Purbhoo for
showing me that one of my early attempts to prove Theorem 10 was
flawed.


\begin{thebibliography}{99}

\bibitem{Chung} F. R. K. Chung, \emph{Pebbling in hypercubes}, SIAM J. Disc. Math \textbf{2} (1989), 467--472.

\bibitem{Glenn1} B. Crull, T. Cundiff, P. Feltman, G. H. Hurlbert,
L. Pudwell, Z. Szaniszlo, Z. Tuza,  \emph{The cover pebbling number
of graphs}, (2004), preprint, http://arxiv.org/abs/math.CO/0406206.

\bibitem{Glenn4} A. Czygrinow, G. Hurlbert, H. A. Kierstead, W. T.
Trotter, \emph{A Note on Graph Pebbling}, Graphs and Combinatorics
\textbf{18} (2002), 219--225.

\bibitem{monthlypaper} A. P. Godbole, N. G. Watson, C. R. Yerger,
\emph{Threshold and complexity results for the cover pebbling game},
(2005), in preparation.

\bibitem{Herv}D. S. Herscovici, \emph{Graham's pebbling conjecture on
products of cycles}, Journal of Graph Theory, \textbf{42} (2003),
141--154.

\bibitem{Glenn2} G. Hurlbert, \emph{A survey of graph pebbling}, Congressus Numerantium \textbf{139} (1999), 41--64.

\bibitem{Glenn3} G. Hurlbert, B. Munyan, \emph{The cover pebbling number of
hypercubes,} (2004), in http://arxiv.org/abs/math.CO/0409368.


\bibitem{Karp} R. M. Karp, {\it Reducibility among combinatorial problems,}
in {\it Complexity of computer computations,} (Proc. Sympos. IBM
Thomas J. Watson Res. Center, Yorktown Heights, N.Y). New Plenum,
New York, pp. 85--103, 1972.

\bibitem{Milans} K. Milans, B. Clark, \emph{The Complexity of Graph
Pebbling}, (2005), preprint, http://arxiv.org/abs/math.CO/0503698.

\bibitem{Vuong} A. Voung, M. I. Wyckoff, \emph{Conditions for
weighted cover pebbling of graphs,} (2004), preprint,
http://arxiv.org/abs/math.CO/0410410.

\bibitem{Watson} N. G. Watson, C. R. Yerger, \emph{Cover pebbling numbers and bounds for certain families of
graphs}, (2004), preprint, http://arxiv.org/abs/math.CO/0409321.

%\bibitem{brown}
%T. C. Brown, P. Erd\H{o}s, and A. R. Freedman, Quasi-progressions
%and descending waves, \emph{J. Combin. Theory Ser. A} \textbf{53},
%no. 1 (1990), 81--95.




\end{thebibliography}
\end{document}